\documentclass{article}

\usepackage{amssymb}

\long\def\rests#1{}

\def\noi{\noindent}

\def\Pf{\noi{\bf Proof.\ \,}}
\def\eop{\hfill\framebox[2.4mm][t1]{\phantom{x}} \vskip 0.15cm } 

\def\vo{vertex operator\ }
\def\voa{vertex operator algebra\ }
\def\voas{vertex operator algebras\ }
\def\svoa{vertex operator superalgebra\ }
\def\Aut{{\rm Aut}}

\newcommand{\sfr}[2]{\leavevmode\kern-.05em
  \raise.5ex\hbox{\the\scriptfont0 #1}\kern-.1em
  /\kern-.15em\lower.25ex\hbox{\the\scriptfont0 #2}\kern.02em}

\def\B{{\mathbb B}}
\def\M{{\mathbb M}}
\def\C{{\cal C}}
\def\D{{\cal D}}
\def\G{{\mathfrak{B}}}
\def\La{L(\sfr{1}{2},0)}
\def\Lb{L(\sfr{1}{2},\sfr{1}{2})}
\def\Lc{L(\sfr{1}{2},\sfr{1}{16})}
\def\VB{V\!B^{\natural}}
\def\V{V^{\natural}}
\def\Wa{W(0)}
\def\Wb{W(\hbox{$\sfr{1}{2}$})}
\def\Wc{W(\hbox{$\sfr{1}{16}$})}

\def\Co{{\rm Co}}

\title{The Group of Symmetries of the 
\hbox{shorter Moonshine Module}}

\author{Gerald~H\"ohn\thanks{
 Department of Mathematics, Kansas State University,
138 Cardwell Hall, Manhattan, KS 66506-2602,
USA. E-mail: {\tt gerald@math.ksu.edu}}}

\date{}

\begin{document}

\bibliographystyle{amsalpha} 

\newtheorem{thm}{Theorem}
\newtheorem{prop}[thm]{Proposition}
\newtheorem{lem}[thm]{Lemma}
\newtheorem{rem}[thm]{Remark}

\newtheorem{cor}[thm]{Corollary}
\newtheorem{conj}[thm]{Conjecture}
\newtheorem{de}[thm]{Definition}
\newtheorem{nota}[thm]{Notation}

\maketitle

\begin{abstract}
It is shown that the automorphism group of the shorter Moonshine module $\VB$
constructed in~\cite{Ho-dr} 
is the direct product of the finite simple group known as the
Baby Monster and the cyclic group of order~$2$. 

\end{abstract}  


\section{Introduction}

The {\it shorter Moonshine module\/} or 
Baby Monster vertex operator superalgebra denoted by $\VB$
was constructed in~\cite{Ho-dr}. It is a 
\svoa of central charge $23\frac{1}{2}$ on which a direct product of 
the cyclic group of order~$2$ and
Fischer's Baby Monster $\B$ acts by automorphisms (\cite{Ho-dr}, Th.~2.4.7).
Furthermore, it was shown loc.~cit.~that the complete automorphism group $\Aut(\VB)$
of $\VB$ is finite. It was conjectured that this is already the full automorphism
group, i.e.~$\Aut(\VB)=2\times \B$.
The purpose of this note is to present a proof for this conjecture.
\begin{thm}\label{main}
The automorphism group of the shorter Moonshine module $\VB$ is
the direct product of the Baby Monster $\B$ and a cyclic group of
order $2$.
\end{thm}
The Baby Monster was discovered by Bernd Fischer in 1973 (unpublished)
and a computer proof of its existence and uniqueness was given
by J.~S.~Leon and C.~C.~Sims in 1977 
(for an announcement see~\cite{LeSi}). 
An independent computer-free construction was given by  R.~Griess
during his construction of the Monster (see~\cite{Gr}).
The Baby Monster is the second largest sporadic simple group 
and has order 
$2^{41}3^{13}5^67^211\cdot 13\cdot 17 \cdot 19 \cdot 23 \cdot 31 \cdot 47$.

\smallskip 

The paper is organized as follows. In the rest of the introduction
we will  recall the definition of $\VB$ from~\cite{Ho-dr} and discuss
the original proof of the corresponding theorem for the Moonshine module:
\begin{thm}\label{mainmonster}
The automorphism group of the Moonshine module $\V$ is the
Monster~$\M$.
\end{thm}
The second section contains the proof of Theorem~\ref{main} which 
is divided into two steps. In the first step, we show that 
$\B$ equals the automorphism group of the algebra $\G'$, 
the subspace of the shorter Moonshine module $\VB$ spanned by
the vectors of conformal weight~$2$. This is done by
identifying certain idempotents in $\G$ and $\G'$ with 
certain conjugacy classes of involutions in $\M$ resp.~$\B$. 
The second step consists in showing that 
the  vertex operator subalgebra $\VB_{(0)}$ 
of even vectors in $\VB$ is generated as 
a \voa by the conformal weight $2$ component~$\G'$. For this, 
we use the theory of framed \voas as developed in~\cite{DGH-virs}.

\bigskip

The Moonshine module $\V$ is
a \voa of central charge~$24$ having
as conformal character
$\chi_{\V}=q^{-1}\sum_{n=0}^{\infty}{\rm dim}\, \V_n\, q^n$ the elliptic
modular function $j-744=q^{-1}+ 196884\, q +  21493760\,q^2 +\cdots $.
It was constructed in~\cite{FLM} (see also~\cite{Bo-ur}) as a
${\bf Z}_2$-orbifold of the lattice \voa $V_{\Lambda}$
associated to the Leech lattice $\Lambda$. 
The algebra  induced on the conformal weight $2$ part $\V_2$ can be 
identified with the Griess algebra $\G$ (see~\cite{FLM}, Prop.~10.3.6). 

We outline the proof of Theorem~\ref{mainmonster}: From the construction
of $\V$ one can find a subgroup of type $2^{1+24}_+.\Co_1$ as well as an extra 
triality involution $\sigma$ inside $\Aut(\V)$ and $\Aut(\G)$
(see~\cite{FLM}, Th.~10.4.12 and Prop.~12.2.9).
It was shown in~\cite{Gr} that $\Aut(\G)$ is finite and certain
theorems from finite group theory (\cite{Smith,GMS}) were invoked to conclude
that the subgroup $\langle 2^{1+24}_+.\Co_1,\,\sigma \rangle$ of $\Aut(\G)$ 
must be the finite simple group called the Monster. 
The finiteness proof for $\Aut(\G)$ was simplified in~\cite{co-monster,Ti-monster}.
Furthermore, Tits showed in~\cite{Ti-monster} that the Monster is the full
automorphism group of $\G$. The proof of Theorem~\ref{mainmonster} is
completed by showing that $\V$ is generated as a \voa by $\G$ 
(see~\cite{FLM}, Th.~12.3.1~(g)). An independent proof of Theorem~\ref{mainmonster}
was given by Miyamoto~\cite{Mia-2004}.

\smallskip

The Monster group has two conjugacy classes of involutions,
denoted by $2A$ and $2B$ in the notation of the Atlas of finite groups
(cf.~\cite{atlas}).
For any involution in  class $2A$ there exists an associated
idempotent $e$ in the Griess algebra~$\G$ called 
the transposition {\it axis\/} of that involution in~\cite{co-monster}.
(In~\cite{co-monster}, $64\,e$ was used, but we prefer to work with the 
idempotents.)
The vertex operator subalgebra of $\V$ generated by $e$ is isomorphic to 
the simple Virasoro \voa $\La$ of central charge $\frac{1}{2}$ 
and has the Virasoro element $2e$ (cf.~\cite{DMZ}, see also~\cite{Ho-dr}, Th.~4.1.2 or~\cite{mia}).

The shorter Moonshine module $\VB$ is a \svoa of central 
charge~$23\frac{1}{2}$ with conformal character
$$\chi_{\VB}=q^{-47/48}(1+ 4371\,q^{3/2} + 96256\, q^2 + 
1143745  \, q^{5/2}  +\cdots ).$$
We recall its definition from~\cite{Ho-dr}, Sec.~3.1 and~4.2:
Let $\La$ be the 
vertex operator subalgebra of $\V$ generated by a transposition axis~$e$.
The Virasoro vertex operator algebra $\La$ has three irreducible 
modules $\La$, $\Lb$ and $\Lc$ 
of conformal weights $0$, $\frac{1}{2}$ and~$\frac{1}{16}$ (cf.~\cite{DMZ}).
Let $\VB_{(0)}$ be the commutant of $\La$ in $\V$.
One has a decomposition
\vspace{-0.4mm}
$$\V\,\,=\,\,\VB_{(0)}\otimes \La\,\, \oplus \,\,
\VB_{(1)}\otimes \Lb \,\,\oplus\,\,
\VB_{(2)}\otimes \Lc$$
 with $\VB_{(0)}$-modules $\VB_{(0)}$,  $\VB_{(1)}$ and $\VB_{(2)}$. 
The shorter Moonshine module  is defined
as $\VB=\VB_{(0)}\oplus \VB_{(1)} $.
A \svoa structure on $\VB$ such that $\VB_{(0)}$ resp.~$\VB_{(1)}$
is the even resp.~odd part can be constructed
from the \voa structure of $\V$ (see~\cite{Ho-dr}, Th.~3.1.3).
The centralizer of the class $2A$ involution associated to the idempotent~$e$
is $2.\B$, a $2$-fold cover of the Baby Monster.
There is an induced action of $2.\B$ on
$\VB$ respecting the \svoa structure such that the  central
involution acts trivially (see~\cite{Ho-dr}, Th.~4.2.6).
An extra automorphism of $\VB$ commuting with the $\B$-action
is the involution which acts as the identity on $\VB_{(0)}$
and by multiplication with $-1$ on $\VB_{(1)}$.

A construction of $\VB$ using the $2A$-twisted
module of $\V$ was explained in~\cite{Ya-twist}.

\medskip

{\it Remarks:\/} The shorter Moonshine module can be understood 
as the \voa analog of 
the {\it shorter Leech lattice\/} (which has $2.\Co_2$ as automorphism group,
where $\Co_2$ is the second Conway group)
and the {\it shorter Golay code\/} (which has $M_{22}{:}2$ as automorphism group,
where $M_{22}$ is a Mathieu group) (cf.~\cite{atlas}).

As a framed \vo superalgebra (for the definition see the next section), %
$\VB$ can in principle be constructed without reference to $\V$.
%
This suggests that a natural way to {\it define\/} the
Baby Monster would be to use the symmetry group of $\VB$.

\medskip

{\it Acknowledgment:\/} The author likes to thank the referees for 
carefully reading the manuscript and the suggestion of the first referee
for simplifying the proof of Proposition~\ref{autunten}.


\section{Proof of Theorem}

Let $V$ be a nonnegatively graded rational \voa having a Virasoro \voa $\La$ as subalgebra.
It was shown in~\cite{mia}, Sect.~4, that in this situation
a nontrivial automorphism of $V$ can often be constructed:
There is the decomposition $V=\Wa\oplus \Wb \oplus \Wc$,
where $W(h)$ is the subspace of $V$ generated by all $\La$-modules
isomorphic to $L(\sfr{1}{2},h)$. If $\Wc\not=0$, one defines
an involution~$\tau$ as the identity on $\Wa\oplus \Wb$
and the multiplication with $-1$ on $\Wc$.
If $\Wc=0$ but $\Wb\not=0$, an involution~$\sigma$
is defined as the identity on $\Wa$ and the 
multiplication with $-1$ on $\Wb$. The nontrivial fusion rules 
of the fusion algebra for $\La$ are
$ \Lb\times \Lb =\La$, $\Lb\times\Lc=\Lc$ and $\Lc\times\Lc=\La+\Lb$ (cf.~\cite{DMZ}).
This is the reason that $\tau$~resp.~$\sigma$ are indeed automorphisms of $V$.
The following result was mentioned in~\cite{mia}, p.~547, but no explicit proof was given.
\begin{lem}\label{axis-monster}
Every idempotent of $\G$ generating a simple Virasoro \voa of central charge 
$\frac{1}{2}$ is a transposition axis. In particular, the map which 
associates to such an idempotent its Miyamoto involution defines a bijection between
the set of such idempotents and the class $2A$ involutions of the Monster.
\end{lem}
\Pf
Let $e$ be an idempotent of $\G$ as in the lemma, let
$\La\subset \V$ be the simple Virasoro \voa of central charge
$\frac{1}{2}$ it generates, and let $\V=\Wa\oplus \Wb \oplus \Wc$
be the decomposition of $\V$ into $\La$-modules as described above.

First we observe that the case $\Wb=\Wc=0$ cannot occur.  
In this case, the 
algebra $\G$ would be the direct orthogonal sum of 
the one dimensional algebra ${\bf C}\cdot e$ and a complementary subalgebra.
Let $e'=x\cdot e+ u$ be another element in the Monster orbit 
of $e$ where $u$ is in the complementary subalgebra.
We get $e'^2=x^2\cdot e + u^2 = x\cdot e+ u =e'$ and
thus either $x=1$ or $x=0$. In the first case, we
would have $\langle e', e'\rangle =\langle e, e\rangle+
\langle u, u\rangle$ and thus $u=0$ and $e'=e$. It
follows that $x=0$ and $\langle e',e\rangle=0$. Thus
the Monster orbit of $e$ would consist of pairwise orthogonal idempotents
generating simple Virasoro \voas of central charge $1/2$. 
Since the central charge of $\V$
is~$24$, there can be at most $48$ such idempotents.
However, no non-trivial factor group of the Monster is a subgroup of any symmetric group
of degree $n\leq 48$.  Hence the Monster fixes $e$ and thus also $1-e$
where $1$ is the identity element of $\G$.
Since the trivial Monster representation occurs only with
multiplicity one, this is impossible. Therefore this case cannot occur
and the idempotent~$e$ defines an involution \hbox{$t\in \Aut(\V)=\M$.}
 
\smallskip
Let $\G=X^+\oplus X^-$ be the decomposition of $\G$ into the $+1$ 
and $-1$ eigenspaces of $t$ and
let $C$ be the centralizer of $t$ in the Monster.
The fusion rules for $\La$ imply that
the algebra multiplication restricts to a $C$-equivariant map
$\mu: X^+\longrightarrow {\rm End}(X^-)$ defined by $\mu(a)x=a_{(1)}x$ for $a\in X^+$ and
$x\in X^-$. 
We claim that \hbox{${\rm Ker}\, \mu=0$.} For,
consider a  $C$-irreducible component of $ {\rm Ker}\, \mu$. 
There are two cases:

\noindent (1) The involution $t$ is in class $2B$ 
and $C\cong 2^{1+24}_+.\Co_1$.
The irreducible $C$-components of $X^+$ have
dimensions $1$, $299$ and $98280$. By the construction of the  
Griess algebra~\cite{Gr}, the $1$-, $299$- and 
$98280$-dimensional component of $X^{+}$
act nontrivially on the $98304$-dimensional component $X^{-}$.
This proves the claim in the $2B$ case.

\noindent (2) The involution $t$ is in class $2A$
and $C\cong 2.\B$.
The decomposition of the Griess algebra $\G$ into irreducible 
components for the centralizer $2.\B$ is shown in
table~\ref{griesszerlegung} (cf.~\cite{Gr,MeNe}).
The components of $X^+$ have dimensions $1$, $1$, $4371$
and $96255$.
Again, it can easily been seen that the $4371$- and the
$96255$-dimensional component act nontrivially
on $X^-$. For example, one can use
the Virasoro frame decomposition of $\V$ given in~\cite{DGH-virs}
(for the notation, see the discussion before Proposition~\ref{auttotal}):
Let $F=\{\omega_1,\ldots,\omega_{48}\}$ be the Virasoro frame for $\V$ as in~\cite{DGH-virs}. 
We can choose $\omega_1=2e$. 
As seen from table~\ref{griesszerlegung}, $\omega_1$ and $\omega_1+\cdots+\omega_{48}$ generate the trivial 
isotypical $2.\B$-module component of $X^-$ and act by multiplication with $1/16$ resp.~$2$ on $X^-$.
The element $\omega_2$ which is contained in the direct sum of the trivial with the $96255$-dimensional
component acts on $X^-$. On the
Viraroso submodule $\Lc\otimes \Lb\otimes (\Lc \otimes \La)^{\otimes 23}$ generated
by a vector $v$ in $X^-$ the action of $\omega_2$ is multiplication with $1/2$, on the Viraroso submodule 
$\Lc\otimes \La\otimes \Lc\otimes  \Lb \otimes \hbox{$ (\Lc \otimes \La)^{\otimes 22}$}$ 
generated by a vector $v'$ in $X^-$ the action of $\omega_2$ is multiplication with~$0$. 
Thus \hbox{$\mu(\omega_2)\not\in {\bf C}\cdot{\rm id}_{X^-}$}
and we see that the action of the $96255$-dimensional component must be nontrivial. 
For the $4371$-dimensional component of $X^+$, choose again a Viraroso submodule 
$ \Lb\otimes \Lb^{\otimes 3}\otimes  \La^{\otimes 44} $ which is generated by a highest weight 
vector $w$ in this component. For $X^-$, choose a Viraroso submodule 
$\Lc\otimes \Lb \otimes ( \Lc\otimes \La)^{\otimes 23}$ which 
is generated by a vector $w'$ in $X^-$. One sees from the fusion rules that the 
algebra product $w\cdot w'$ is nontrivial. This finishes the proof of the claim 
in the $2A$ case.

\smallskip
Let  $U=C\cdot e$ be the $C$-submodule of $X^+$ generated by $e$.
Since $e$ acts by multiplication with a non-zero scalar $\lambda$ on $X^+$
(one has $\lambda=\frac{1}{32}$ if $t$ is of type $\tau$ and 
$\lambda=\frac{1}{4}$ if $t$ is of type $\sigma$),
all elements of $U$ are mapped by $\mu$
to ${\bf C}\cdot{\rm id}_{X^-}\subset {\rm End}(X^-)$,
and $\mu(e)=\lambda \cdot {\rm id}_{X^-}$ with $\lambda\not=0$. 
It follows that $U$ must be a trivial representation of~$C$.
In case~(1), this implies that $e$ must be the identity of $\G$.
Since the identity  corresponds to the Virasoro element of central 
charge~$24$, this case is impossible. In case~(2),
we see that~$e$ is contained in the $2$-dimensional trivial
$C$-component ${\bf  C}\oplus {\bf C}$, and this implies that $e$ 
must be an axis since this component is spanned by $e$ and the identity
\nopagebreak
 of~$\G$. \eop
\begin{table}\caption{Decomposition of the Griess algebra $\G$ under $2.\B$}\label{griesszerlegung}
{\small
$$\begin{array}{l|cccccccccc}
\hbox{$2.\B$ representations:} & \underline{1} & \oplus & \underline{1} &
 \oplus & \underline{96255} & \oplus & \underline{4371} & \oplus & 
\underline{96256}  \\ \hline
\hbox{Eigenvalues of $\mu(2e)$:\phantom{\Huge{g$\!\!\!$}}}
 & 2 & & 0 & & 0 & & \frac{1}{2} & & \frac{1}{16} \\ \hline 
\hbox{Eigenvalues of $\tau$:} & +1 & & +1 & & +1 & & +1 & & -1 \\ \hline
 &  & & \multicolumn{3}{c}{\underbrace{\phantom{ ------ }}} & & & & &\\  
\hbox{Parts of $\VB$:} &  & &  \multicolumn{3}{c}{\VB_2=\G'}  & & \VB_{3/2} & & \\ 
\end{array}$$
}
\end{table}

Denote with $\G'$ the subspace of $\VB$ 
spanned by the vectors of conformal weight~$2$.
The vertex operator algebra structure on $\VB_{(0)}$ induces
an algebra structure on $\G'$.
\begin{lem}\label{axis-baby}
The map which associates to an idempotent of $\G'$ generating 
a simple Virasoro \voa of central charge $\frac{1}{2}$
its Miyamoto involution defines a bijection between
the set of transposition axis contained in $\G'$
and the class $2B$ involutions of the Baby Monster.
\end{lem}
\Pf Since, by construction, $\G'$ is a subalgebra of $\G$, 
every idempotent of $\G'$ is also an idempotent of $\G$. 
By Lemma~\ref{axis-monster}, the idempotents which generate 
a simple Virasoro \voa of central charge $\frac{1}{2}$ must 
be transposition axes.

The set of $2A$ involutions of $\M$ decomposes under the 
conjugation action of the centralizer $2.\B$ of a fixed 
$2A$-involution~$\tau$
into nine orbits~\cite{No,GMS}. These orbits can be labeled by the 
Monster conjugacy class of the product of an involution 
in an orbit with the fixed $2A$-involution $\tau$ corresponding to
an axis~$e$.
The orbit for the class $2B$ corresponds to the axis of $\G$
contained in the $\G'$ component of $\G$ (see~\cite{MeNe}, Th.~5).
The stabilizer in $\B$ of an involution in this orbit  
(the central element $2$ in $2.\B$ acts trivial) is $2^{1+22}_+.\Co_2$ 
(cf.~\cite{GMS}) and the action of $\B$ on the orbit
can be identified with the conjugation action of $\B$  
on the class of $2B$ involutions of $\B$.

Using the map
$2.\B={\rm Cent}(\tau)={\rm Stab}_\M(e)\longrightarrow {\rm Aut}(\VB_{(0)})$,
$g\mapsto g|_{{\rm Aut}(\VB_{(0)})}$,
with kernel $\langle \tau  \rangle$, we see that 
the involution $\tau'$ on $\VB_{(0)}$ associated to an idempotent 
$e$ in $\G'\subset \G$ is the restriction to $\VB_{(0)}$
of the associated involution $\tau$ of~$\V$. \phantom{xxxxxxxxxx}\hfill
\eop

\begin{prop}\label{autunten}
The automorphism group of the algebra $\G'$ is the Baby Monster.
\end{prop}
\Pf
Let $G=\Aut(\G')\geq\B$. Since all the idempotents $e$ of $\G'$ generating a
central charge $\frac{1}{2}$ Virasoro algebra correspond to the class 
$2B$ involutions of $\B$, the conjugation action of $G$ 
leaves the set of $2B$ involutions invariant. 
Since $\B$ is generated by its $2B$ involutions it is a normal subgroup of $G$. 
Since ${\rm Aut}(\B)=\B$ (cf.~\cite{atlas}), we see $G= A.\B$  and
$A$ centralizes $\B$. Recall from Table~\ref{griesszerlegung} that 
$\G'=\underline{{\bf 1}} \oplus \underline{{\bf 96255}}$ as $\B$-module.
Since $G$ fixes the Virasoro element  $\omega$  of $\VB_{(0)}$ this is
also a decomposition of $G$-modules.
Let $a\in A$. Then $a$ must act by a scalar $\lambda$ on the 
$96255$-dimensional component. We write an idempotent~$e$
as above in the form $e=x \cdot \omega+ u$ 
where $\omega$ is the Virasoro element and $u$ is in the $96255$-dimensional component. 
Clearly, $u\not=0$. We obtain
$a(e)=x\cdot \omega+ \lambda u$ since $a$ fixes $\omega$.
As $e=a(e)$ it follows that $\lambda=1$. Thus $a=1$ and $A$ is trivial.
\phantom{xxxxxxxx} \eop


\smallskip

{\it Remark:\/} In a similar way, one could prove that $\Aut(\G)=\M$.
However, the proof of Lemma~\ref{axis-monster} already assumes that 
the Monster is the whole automorphism group of $\Aut(\G)$.

\medskip

For proving that $\VB_{(0)}$ is generated by the subspace $\G'$, 
we recall some results about framed \voas from~\cite{DGH-virs}.
A subset $F=\{\omega_1,\ldots,\omega_r\}$ of a simple vertex operator algebra
$V$ is called a {\em Virasoro frame\/}  if the $\omega_i$
for \hbox{$i=1$, $\ldots$, $r$} generate mutually commuting simple Virasoro
vertex operator algebras  $\La$ of central 
charge $1/2$~and $\omega_1+\cdots +\omega_r$ is the Virasoro element of~$V$. 
Such a \voa $V$ is called a {\em framed vertex operator algebra\/}.
Under $\La^{\otimes r}$, the vertex operator subalgebra 
spanned by the Virasoro frame, $V$ decomposes as
$$V=\bigoplus_{h_1,\,\ldots,\,h_r\in\{0,\,\sfr{1}{2},\, \sfr{1}{16}\}} 
n_{h_1,\,\ldots,\,h_r}\, L(h_1,\,\ldots,\,h_r),$$
into isotypical components of modules
\hbox{$ L(h_1,\,\ldots,\,h_r)=L(\sfr{1}{2},h_1)\otimes\cdots \otimes L(\sfr{1}{2},h_r)$.}
Using this decomposition, one defines the two binary codes 
$${\cal C}=\{c\in{\bf F}_2^r \mid n_{c_1/2,\,\ldots,\,c_r/2}\not=0\}\ 
\hbox{and}\   
{\cal D}=\{d\in{\bf F}_2^r \mid 
\sum_{ h_1,\ldots,\, h_r}  n_{h_1,\,\ldots,\,h_r}\not=0 \},$$
where the sum in the definition of ${\cal D}$ runs over those $h_i$ with $h_i\in \{0,\, \sfr{1}{2}\}$ if $d_i=0$ 
and $h_i=\sfr{1}{16}$ if $d_i=1$.
It was shown in~\cite{DGH-virs} that both codes
are linear, ${\cal D}\subset {\cal C}^{\perp}$ ($C^{\perp}$ denotes the code orthogonal 
to $C$), ${\cal C}$ is even,
and all weights of vectors in ${\cal D}$ are divisible by~$8$.
For $c\in {\cal C}$, we denote the sum of all irreducible modules
isomorphic to $L(c_1/2,\,\ldots,\,c_r/2)$ by $V(c)$; for 
$I\in {\cal D}$, we denote the sum of all irreducible modules
isomorphic to  $L(h_1,\,\ldots,\,h_r)$ such that
$h_i\in \{0,\, \sfr{1}{2}\}$ if $I_i=0$  and
$h_i=\sfr{1}{16}$ if $I_i=1$ for \hbox{$i=1$, $\ldots$, $r$} by $V^I$. 
\begin{prop}\label{auttotal}
The even part $\VB_{(0)}$ of the shorter Moonshine module $\VB$ 
is generated by  $\G'$, the subspace of $\VB$
consisting vectors of conformal weight~$2$.
\end{prop}
\Pf The Moonshine module $\V$ is a framed \voa
of central charge~$24$ (cf.~\cite{DMZ}). To any Virasoro frame in $\V$ there are
the two associated binary codes $\C$ and $\D$ of length~$48$. 
For the construction of the shorter Moonshine module,
one has to fix a simple central charge~$\frac{1}{2}$ Virasoro
vertex operator subalgebra generated by an axis.
It follows from Lemma~\ref{axis-monster} that
all such vertex operator subalgebras 
are equivalent, one can choose an arbitrary $\La$ of a given
Virasoro frame in $\V$ and the remaining $47$ factors $\La$
define a Virasoro frame of~$\VB_{(0)}$.
The associated binary codes $\C'$ and $\D'$ for 
the framed \voa~$\VB_{(0)}$ are obtained from the codes
$\C$ resp.~$\D$ by taking those vectors 
of $\C$ resp.~$\D$ which have the value $0$ at the
coordinate corresponding to the first chosen $\La$
and then omitting this coordinate.

We fix now the Virasoro frame of $\V$ studied in~\cite{DGH-virs}, Sect.~5, 
Example~II.
The code $\D$ has in this case the generator matrix 
$$\left(\begin{array}{ccc}
1111 1111 1111 1111 & 0000 0000 0000 0000 & 0000 0000 0000 0000 \\
0000 0000 0000 0000 & 1111 1111 1111 1111 & 0000 0000 0000 0000 \\
0000 0000 0000 0000 & 0000 0000 0000 0000 & 1111 1111 1111 1111 \\ 
0000 0000 1111 1111 & 0000 0000 1111 1111 & 0000 0000 1111 1111 \\
0000 1111 0000 1111 & 0000 1111 0000 1111 & 0000 1111 0000 1111 \\
0011 0011 0011 0011 & 0011 0011 0011 0011 & 0011 0011 0011 0011 \\
0101 0101 0101 0101 & 0101 0101 0101 0101 & 0101 0101 0101 0101 
\end{array}\right)_{\textstyle }$$
and $\C$ is the dual code of $\D$, i.e.,~$\C=\{c\in {\bf F}_2^{48}\mid \sum_{i=1}^{48}
c_i\cdot d_i=0\ \hbox{for all}\  d\in \D \}$ (Theorem~5.8 of~\cite{DGH-virs}). 
In particular, $\C$ is even.
It was shown in~\cite{DGH-virs} that $\C$ has minimal weight~$4$.

We claim that the codewords of weight~$4$ in $\C$ span the code:
The $48$ coordinates of $\C$ are naturally partitioned into three blocks
of $16$ coordinates by the first three rows of the above matrix. 
Let $c$ be a codeword of $\C$ of weight~$6$ or larger. By the pigeon hole 
principle, we can find two blocks such that the support of~$c$ has at least four 
coordinates at the $32$ coordinates of these two blocks. The codewords of
$\C$ having the value $0$ at the remaining third block can be identified
with the extended Hamming code ${\cal H}_{32}$ of length~$32$ since its dual code 
is the extended simplex code of that length whose generator matrix equals 
the restriction of $\D$ to the two blocks. 
The code ${\cal H}_{32}$ has the property that for any set of three of its $32$ 
coordinates there is a codeword~$v$ of weight~$4$ 
in ${\cal H}_{32}\subset \C$ containing this set of
coordinates in its support, i.e., the codewords of weight~$4$ define a
Steiner system $S(3,4,32)$ (cf.~\cite{McS}, Ch.~2, Thm.~15).
By choosing the three coordinates from
the support of $c$ contained in the  two selected  blocks, we see
that $c+v$ is a codeword of $\C$ of smaller weight then the weight of~$c$. 
This implies the claim. 

Now we claim that the codewords of weight~$4$ in $\C'$ span $\C'$:  
By using the claim for ${\cal C}$, we see that a codeword $c$ of $\C$ with
value~$0$ at a fixed coordinate is the sum of codewords in $\C$ of 
weight~$4$ such that in the sum an even number of the codewords have the value
$1$ at the fixed coordinate. The sum $d$ of two different weight~$4$ codewords having the 
value $1$ at the fixed coordinate is a codeword of weight $4$ or $6$.
We must show that $d$ is the sum of weight~$4$ codewords of $\C'$. So 
if $d$ has weight~$4$ we are done.  If $d$ has weight $6$,
let $i\in\{0,\,1,\,\ldots,\,6\}$ 
be the number of coordinates of the support of $d$ 
which are in the same block as the fixed
coordinate used in the definition of $\C'$. 

For $i\leq 3$,
the support of $d$ contains $6-i\geq 3$ coordinates in the two other blocks.
In this case we can choose a weight~$4$ vector $v$ as above with support in these two blocks
and obtain with $v$ and $d+v$ two weight $4$ vectors in $\C'$ with sum $d$.

For $i\geq 4$, one observes from the definition of $\C$ that the codewords of
$\C$ having the value $0$ at the two other block 
can be identified with the extended Hamming code ${\cal H}_{16}$ of length~16:
its dual code is the extended simplex code of that length whose generator matrix equals the
restriction of $\D$ to the remaining block. 
We choose now in two different ways $3$ coordinates out of a set of $4$ coordinates
contained in the $i\geq 4$ coordinates which form the intersection
of the support of~$d$ with the block containing the distinguished coordinate.
Since the weight~$4$ codewords of ${\cal H}_{16}$ form a Steiner system $S(3,4,16)$, 
there are two codewords $v$ and $v'$ of weight~$4$ in ${\cal H}_{16}$
containing these two sets of $3$ coordinates in their support. Not both of them can have 
the value~$1$ at the distinguished coordinate since otherwise their sum would have weight~$2$.
So we have found with $v$ and $d+v$ or with $v'$ and $d+v'$ two vectors of 
weight~$4$  in $\C'$ with sum $d$ and we are done again.

It is clear that $\D$ and $\D'$ are 
both generated by their codewords of weight~$16$ and~$24$.

\smallskip 
In~\cite{DGH-virs}, Prop.~2.5.~(5) and~(6), it was shown that
for a framed \voa $V$
the module $V(c)$  generated by the isotypical
component corresponding to a codeword $c$
in the associated binary code~$\C$ is irreducible and 
the span of the images of the maps from $V(c)\times V(c')$
to $V(c+c')$, $c$, $c'\in \C$, defined  by the components of the
vertex operator generate $V(c+c')$.
Since the Virasoro frame is contained in $V_2$
and the conformal weight of a Virasoro highest weight
vector in $V(c)$ for $c\in \C$ is half the weight of $c$, it follows 
from the above discussion that the \vo subalgebra 
$V^0=\bigoplus_{c\in \C}V(c)$ for $V$ equal to $\V$ resp.~$\VB_{(0)}$ 
is generated by $\G$ resp.~$\G'$.
 
Similarly, using~\cite{DGH-virs}, Prop.~2.5.~(1) and (2), one has:
$V=\bigoplus_{I\in \D}V^I$, the $V^I$ are irreducible $V^0$-modules,
and the span of the images of the maps from $V^I\times V^{I'}$
to $V^{I+I'}$, $I$, $I'\in \D$, defined by the components of the 
vertex operator generate $V^{I+I'}$. By direct inspection of the
Virasoro decomposition of the moonshine module given explicitly
in~\cite{DGH-virs}, one finds for any codeword~$I$ of weight~$16$ or~$24$ in
$\D$ a Virasoro highest weight vector in $(\V)^I_2$. For example,
there is the Virasoro highest weight module 
$L(\sfr{1}{2},0^{15},\sfr{1}{2},0^{15},(\sfr{1}{16})^{16})$. (In more detail,
one can use the decomposition polynomial given in~\cite{DGH-virs}, Cor.~5.9
and observe that the situation is symmetric for all codewords in $\D$ of
fixed weight.) For $\VB_{(0)}$, one observes that for all $I\in \D'$
the Virasoro highest weight module can be chosen to have 
conformal weight~$0$ at the fixed coordinate. (Again,
this can be seen from the decomposition polynomial.)

Putting the above arguments together, we have proven that $\VB_{(0)}$ 
is generated by~$\G'$.
\eop

\smallskip

{\it Remark:\/} The proof of Prop.~\ref{auttotal} also shows that
the Moonshine module $\V$ is generated by the Griess algebra $\G$,
a result first proven in~\cite{FLM}, Prop.~12.3.1~(g) by using similar 
ideas.

\medskip

Theorem~\ref{main} follows now from Prop.~\ref{autunten} 
and Prop.~\ref{auttotal} as explained in~\cite{Ho-dr}, Th.~4.2.7:
Any automorphism $g$ of $\VB$ must fix the even and odd components 
$\VB_{(0)}$ and $\VB_{(1)}$, i.e., can be written in the form 
$g=g_0\oplus g_1$ and one has $g_0\in\Aut(\VB_{(0)}) = \B$. 
Let $g'=g_0'\oplus g_1'$ be another automorphism with $g_0=g_0'$.
Since $\VB_{(1)}$ is an irreducible $\VB_{(0)}$-module
(\cite{Ho-dr}, Th.~4.2.2) and $g^{-1}g'={\rm id}\oplus g_1^{-1}g_1'$
is a $\VB_{(0)}$-module homomorphism, there is a scalar~$s$ such that 
$g'=g_0\oplus s\cdot g_1$. Since $\VB$ is a vertex operator 
superalgebra, one gets $s=\pm 1$. \eop

\medskip



\small


\newcommand{\etalchar}[1]{$^{#1}$}
\providecommand{\bysame}{\leavevmode\hbox to3em{\hrulefill}\thinspace}
\providecommand{\MR}{\relax\ifhmode\unskip\space\fi MR }
\providecommand{\MRhref}[2]{%
  \href{http://www.ams.org/mathscinet-getitem?mr=#1}{#2}
}
\providecommand{\href}[2]{#2}

\end{document}